\documentclass[a4paper,12pt,reqno]{amsart}
\usepackage{graphics}
\newtheorem{lem}{Lemma}

\newtheorem{theo}{Theorem}

\numberwithin{equation}{section}

\newcommand{\id}{\operatorname{id}}

\newcommand{\pbinom}[2]{\left[ #1 \atop #2 \right]_P}
\newcommand{\pe}{{\,_P E}}
\newcommand{\pd}{{\, _P \Delta}}
\newcommand{\pqe}{{\,_P ^Q E}}
\newcommand{\pqid}{{\,\,_P ^Q \id}}

\begin{document}

\title[Operator formula]{The operator formula for monotone triangles -- 
simplified proof and three generalizations}

\author[Ilse Fischer]{Ilse Fischer}

\thanks{Supported by the Austrian Science Foundation
    (FWF), grant number S9607--N13, in the framework of the National Research
    Network  ``Analytic Combinatorics and Probablistic Number Theory''.}

\begin{abstract}
We provide a simplified proof of our operator formula for the number of monotone triangles 
with prescribed bottom row, which enables us to deduce three generalizations
of the formula. One of the generalizations concerns a certain weighted enumeration of 
monotone triangles which specializes to the weighted enumeration of
alternating sign matrices with respect to the number of $-1$s in the matrix
when prescribing $(1,2,\ldots,n)$ as the bottom row of the monotone triangle.
\end{abstract}

\maketitle

\section{Introduction}

A {\it monotone triangle} is a triangular array of integers of the following form
\begin{center}
\begin{tabular}{ccccccccccccccccc}
  &   &   &   &   &   &   &   & $a_{n,n}$ &   &   &   &   &   &   &   & \\
  &   &   &   &   &   &   & $a_{n-1,n-1}$ &   & $a_{n-1,n}$ &   &   &   &   &   &   & \\
  &   &   &   &   &   & $\dots$ &   & $\dots$ &   & $\dots$ &   &   &   &   &   & \\
  &   &   &   &   & $a_{3,3}$ &   & $\dots$ &   & $\dots$ &   & $a_{3,n}$ &   &   &   &   & \\
  &   &   &   & $a_{2,2}$ &   & $a_{2,3}$ &  &   $\dots$ &   & $\dots$   &  & $a_{2,n}$  &   &   &   & \\
  &   &   & $a_{1,1}$ &   & $a_{1,2}$ &   & $a_{1,3}$ &   & $\dots$ &   & $\dots$ &   & $a_{1,n}$ &   &   &
\end{tabular},
\end{center}
which is monotone increasing in northeast and in southeast direction and strictly increasing along rows, that is $a_{i,j} \le a_{i+1,j+1}$ for $1 \le i \le j < n$,
$a_{i,j} \le a_{i-1,j}$ for $1 < i \le j \le n$ and $a_{i,j} < a_{i,j+1}$ for $1 \le i \le j \le n-1$. Monotone triangles with bottom row $(1,2,\ldots, n)$ correspond to $n \times n$ alternating sign matrices, the enumeration of which (there are exactly $\prod\limits_{j=0}^{n-1} \frac{(3j+1)!}{(n+j)!}$ of them) provided an open problem for quite some time, see \cite{bressoud}.

\medskip

In \cite{fischermonotone} we have shown that the number of monotone 
triangles with bottom row $(k_1,\ldots,k_n)$ is given by
\begin{equation}
\label{formula}
\prod_{1 \le s < t \le n} (\id - E_{k_s} +  E_{k_s} E_{k_t}) 
\prod_{1 \le i < j \le n} \frac{k_j - k_i}{j-i},
\end{equation}
where $E_x$ denotes the shift operator, defined as $E_x p(x) = p(x+1)$. 
In such an ``operator formula'', the product of operators is understood as the composition. 
Moreover note that the shift operators with respect to different variables commute and, 
consequently, it does not matter in which order the operators in the product 
are applied. This formula was the basis for our new proof of the refined alternating sign matrix theorem \cite{fischernewproof}.

\medskip

The purpose of this article is to present a (very much) simplified proof of \eqref{formula} and three
generalizations that arise quite naturally out of this new proof. Next we describe these generalizations.

\medskip

Regarding the first generalization consider the following inverse question. Given a polynomial $r(X,Y)$ in $X, X^{-1}, Y, Y^{-1}$ (e.g., over $\mathbb{C}$), find a combinatorial interpretation for the numbers
\begin{equation}
\label{generalform}
\prod_{1 \le s < t \le n} r(E_{k_s}, E_{k_t}) 
\prod_{1 \le i < j \le n} \frac{k_j - k_i}{j-i}.
\end{equation}
For the polynomial $r(X,Y)=\id- X + X Y$, a combinatorial interpretation is obviously given by monotone triangles with prescribed bottom row; for the polynomial $r(X,Y)= Y$, a combinatorial interpretation is given by Gelfand--Tsetlin patterns with bottom row 
$(k_1,\ldots,k_n)$ (which are almost defined as monotone triangles only they need not 
necessarily strictly increase along rows) as the number of these objects is given by
\begin{equation}
\label{semistandard}
\prod_{1 \le i < j \le n} \frac{k_j - k_i + j-i}{j - i}, 
\end{equation}
see, for instance, \cite{fischermethod}. Note that Gelfand--Tsetlin patterns with bottom row $(k_1,\ldots,k_n)$ are 
in bijection with semistandard tableaux of shape $(k_n,k_{n-1},\ldots,k_1)$. Here, we present a combinatorial interpretation of \eqref{generalform} for all polynomials $r(X,Y)$ of the form 
\begin{equation}
\label{specialform}
r(X,Y) = Y + (X-\id) \, (Y-\id) \, s(X,Y),
\end{equation}
where $s(X,Y)$ is an arbitrary polynomial in $X, X^{-1}, Y, Y^{-1}$. (Obviously, \eqref{formula} is the special case 
$s(X,Y) = 1$ and \eqref{semistandard} is the special case $s(X,Y)=0$.) Notably, this interpretation is not only valid for
$(k_1,\ldots,k_n) \in \mathbb{Z}^n$ with $k_1 < k_2 < \ldots < k_n$, but also
for general $(k_1,\ldots,k_n) \in \mathbb{Z}^{n}$ and thus this will also lead to an interpretation of \eqref{formula} for all $(k_1,\ldots,k_n) \in \mathbb{Z}^n$.

\medskip

The second generalization concerns the weighted enumeration of alternating
sign matrices with respect to the number of $-1$s in the matrix, which was
introduced by Mills, Robbins and Rumsey, see  \cite{mills}. A $-1$ in the
alternating sign matrix translates into an entry $a_{i,j}$ with $i \not= 1$
and $a_{i-1,j-1} < a_{i,j} < a_{i-1,j}$ in the corresponding monotone triangle
and, therefore, we define the $Q$--weight of a monotone triangle as $Q$
raised to
the power of the number of such entries. 

\begin{theo}
\label{Qweight}
The generating function of monotone triangles with prescribed bottom row
$(k_1,\ldots,k_n)$, $k_{1} < k_{2} < \cdots < k_{n}$, and with respect to the 
$Q$--weight is given by 
$$
\prod_{1 \le s < t \le n} \left( \id - (2-Q) E_{k_s} + E_{k_s} E_{k_t}  \right) 
\prod_{1 \le i < j \le n} \frac{k_j - k_i}{j - i}.
$$
\end{theo}

Notably, this generating function has already been introduced by Mills, Robbins and Rumsey in \cite[Section 5]{mills}.
(However, Theorem~\ref{Qweight} provides the first formula for this generating function.)

\medskip

Finally, the third generalization is a weighted enumeration of a certain type
of Gelfand--Tsetlin patterns (which we denote as {\it weak monotone triangles} since
they include monotone triangles) with prescribed bottom 
row, which reduces to the enumeration of monotone triangles if we take the limit $P \to 1$. (That is for $P=1$ the 
weight is $1$ for monotone triangles and $0$ for weak monotone triangles that are not strictly increasing along 
all rows.) A {\it weak monotone triangle} is a Gelfand--Tsetlin pattern $(a_{i,j})_{1
  \le i \le j \le n}$ with $a_{i,j-1} < a_{i-1,j}$ if $i \not= 1$ and $i <
j$. The respective weight is defined as 
$$
\prod_{a_{i,j} : a_{i,j} < a_{i-1,j}} (P^{a_{i,j}} - [a_{i,j} = a_{i,j-1}]),
$$
where $[\text{statement}]=1$ if the statement is true and $0$ otherwise.

\medskip 

In order to state this generating function, we need to introduce the following $P$--generalizations of the difference operator $\Delta_x:= E_x - \id$ and the shift operator: the $P$--difference operator is defined as $\pd_x = P^{-x} \Delta_x$ and the $P$--shift operator is defined as $\pe_x = \pd_x + \id$. If we set $P=1$ then we obtain the 
ordinary operators. Note that these operators commute, i.e. $\pd_x \pd_y = \pd_y \pd_x$, 
$\pe_x \pe_y = \pe_y \pe_x$ and $\pe_x \pd_y = \pd_y \pe_x$.

\medskip

\begin{theo} 
\label{Pweight}
The generating function of weak monotone triangles with prescribed bottom row
$(k_1,\ldots,k_n)$, $k_{1} \le k_{2} \le \cdots \le k_{n}$,
and with respect to the $P$--weight is
$$
P^{\binom{n+1}{3}}  \prod_{1 \le s < t \le n} \left( {\pe_{k_t}}  + {\pd_{k_s}} {\pd_{k_t}}  \right) 
\prod_{1 \le i < j \le n} \frac{P^{k_j} - P^{k_i}}{P^j - P^i}.
$$
\end{theo}

The paper is organized as follows. In the following section we introduce a general recursion and a master theorem, 
which implies the three generalizations all at once. In Section~\ref{interpret} we deal with the first generalization and provide 
the combinatorial interpretation for \eqref{generalform} if $r(X,Y)$ is of the form given in \eqref{specialform}. In 
Section~\ref{PQ} we deduce the $Q$--enumeration of monotone triangles and the
$P$--enumeration of weak monotone triangles from the main theorem. Moreover, we derive some results (old and new) for the special case $Q=2$ in 
this section. 
In Section~\ref{proof} we finally prove the master theorem. In Section~\ref{future} we discuss further projects along these lines.

\section{The recursion and the master theorem}
\label{master}

We define a $PQ$--shift operator 
as $\pqe_x = Q \id + \pd_x$ and a $PQ$--identity as $\pqid_{x} = Q \pe_x - \pd_x$. For $Q=1$ we have 
$\pqe_x = \pe_x$ and $\pqid_x = \id_x$ and again these operators commute with each other and also 
with the operators that we have introduced before the statement of Theorem~\ref{Pweight}.
Moreover note that $\pqe_x p(x) = \pqid_{x} p(x) = Q p(x)$ if $p(x)$ is
constant (with respect to $x$).

\medskip

Let $S$ be a finite subset of $\mathbb{Z}^2$ and $f: S \to \mathbb{C}$ be a function. For a given 
$(k_1,\ldots,k_n) \in \mathbb{Z}^n$ and a function $A(l_1,\ldots,l_{n-1})$ on $\mathbb{Z}^{n-1}$, we define the summation operator 
$$
 \sum_{(l_1,\ldots,l_{n-1})}^{(k_1,\ldots,k_n)} A(l_1,\ldots,l_{n-1})
$$
associated to the pair $(S,f)$ by induction with respect to $n$. If $n=0$ then the application of the 
operator gives zero, for $n=1$ we set $\sum\limits_{(-)}^{k_1} A = A$. If $n \ge 2$ then we define 
\begin{multline*}
\sum_{(l_1,\ldots, l_{n-1})}^{(k_1,\ldots,k_n)} A(l_1,\ldots,l_{n-1})  
=  \left( Q^{-1} \pqe_{k_n} \pqid_{k^{*}_{n-1}} 
\sum_{l_{n-1}=k^{*}_{n-1}}^{k_n-1} P^{l_{n-1}} \sum_{(l_1,\ldots,l_{n-2})}^{(k_1,\ldots,k_{n-1})} 
A(l_1,\ldots,l_{n-1}) \right.   \\ \left. + \, \sum_{(l_1,\ldots,l_{n-3})}^{(k_1,\ldots,k_{n-2})}
\sum_{(i,j) \in S} f(i,j) \left. {\pe_{k_{n-1}}}^i {\pe_{k^{*}_{n-1}}}^j A(l_1,\ldots,l_{n-3},k_{n-1},k^{*}_{n-1}) \right)
\right|_{k^{*}_{n-1}=k_{n-1}},
\end{multline*} 
where here and in the 
following  $\sum\limits_{x=a}^{b} g(x): = - \sum\limits_{x=b+1}^{a-1} g(x)$ if $a > b$\footnote{Note that this 
implies $\sum\limits_{x=a}^{a-1} g(x)=0$.} and, for $ i < 0$,  
\begin{equation}
\label{inverse}
 {\pe^{i}_x} = (\id +  \pd _x)^{i} := \left( \sum_{j=0}^{\infty} (-1)^j
   {\pd^j_x} \right)^{-i}.
\end{equation}
(If $i<0$ then we will apply the operator ${\pe^{i}_x}$ only to functions $f(x)$ with  ${\pd^j_x} f(x) = 0$ for a 
$j \ge 0$, i.e. to functions for which the sum in \eqref{inverse} is in fact
finite. Note that the operator 
$\pe_x^{i}$ specializes to $E^{i}_x$ for $P=1$ also if $i<0$.) This
generalizes the summation operator from \cite{fischermonotone}, where we have considered the special case  $P=1, Q=1$, $S=\{(0,0)\}$ and $f(0,0)=-1$. 

\medskip

The $P$--binomial coefficient is defined as 
$$
\pbinom{x}{m} = \frac{(1-P^x)(1-P^{x-1}) \cdots (1-P^{x-m+1})}{(1-P^m)(1-P^{m-1}) \cdots (1-P)}
$$ 
and again we obtain the ordinary binomial coefficient $\binom{x}{m}$ if we take the limit $P \to 1$.
We define $\alpha_{P,Q}(n,m,S,f;k_1,\ldots,k_n)$ inductively with respect to $n$: let 
$\alpha_{P,Q}(1,m,S,f;k_1)=\pbinom{k_1}{m}$ and 
$$
\alpha_{P,Q}(n,m,S,f;k_1,\ldots,k_n) = \sum_{(l_1,\ldots, l_{n-1})}^{(k_1,\ldots,k_n)} \alpha_{P,Q}(n-1,m,S,f;l_1,\ldots,l_{n-1}).
$$
In order to see that $\alpha_{P,Q}(n,m,S,f;k_1,\ldots,k_n)$ is well--defined even if 
$S \not\subseteq \mathbb{Z}_{\ge 0}^2$ (i.e. the application of the operator $\pe_{k_j}^i$ to $\alpha_{P,Q}(n,m,S,f;k_1,\ldots,k_n)$
makes sense even if $i < 0$)
observe that 
\begin{equation}
\label{diff}
{\pd_x} \pbinom{x}{m} = \pbinom{x}{m-1} P^{-m+1}
\end{equation}
and 
\begin{equation}
\label{qsum}
\sum_{x=a}^{b} P^x \pbinom{x}{m} =  \sum_{x=a}^{b} P^{x+m} \pd_x \pbinom{x}{m+1} = 
P^m \left( \pbinom{b+1}{m+1} - \pbinom{a}{m+1} \right)
\end{equation}
implies by induction with respect to $m$ that 
$\alpha_{P,Q}(n,m,S,f;k_1,\ldots,k_n)$ is a linear combination of expressions of the form
$ \pbinom{k_1}{r_1} \pbinom{k_2}{r_2} \cdots \pbinom{k_n}{r_n}$ over $\mathbb{C}[P,P^{-1},Q,Q^{-1}]$, and consequently  
(also by  \eqref{diff}), there is an $i \ge 0$ such that $\pd^i_{k_j} \alpha_{P,Q}(n,m,S,f;k_1,\ldots,k_n)$
vanishes. 

\medskip

Clearly, $\alpha_{1,1}(n,0,\{(0,0)\},-1;k_1,\ldots,k_n)$ is equal to \eqref{formula}. Moreover, it is 
easy to see that $\alpha_{1,1}(n,0,\emptyset,-;k_1,\ldots,k_n)$ is the number of Gelfand--Tsetlin patterns with bottom row 
$(k_1,\ldots,k_n)$ and therefore equal to \eqref{semistandard}. For the general situation we will prove the following theorem.

\begin{theo} 
\label{main}
Let $n$ be a positive integer,  $m$ be a non--negative integer, $S \subseteq \mathbb{Z}^2$ be a finite set and $f: S \to \mathbb{C}$ be a function. Then $\alpha_{P,Q}(n,m,S,f;k_1,\ldots,k_n)$
is given by 
\begin{multline*}
P^{\frac{1}{6} (n-1)(n^2 - 2n + 6m)} Q^{-\binom{n}{2}} \prod_{1 \le s < t \le n} \left( {\pqe_{k_t}} \pqid_{k_s} - {\pd_{k_s}} {\pd_{k_t}} 
\sum_{(i,j) \in S} f(i,j) {\pe_{k_t}}^i {\pe_{k_s}}^j \right) \\ \det_{1 \le i,j \le n}
\pbinom{k_i}{j-1+\delta_{j,n} m}.
\end{multline*}
For $m=0$ this simplifies to
$$
P^{\binom{n+1}{3}}  Q^{-\binom{n}{2}} \prod_{1 \le s < t \le n} \left( {\pqe_{k_t}} \pqid_{k_s} - {\pd_{k_s}} {\pd_{k_t}} 
\sum_{(i,j) \in S} f(i,j) {\pe_{k_t}}^i {\pe_{k_s}}^j \right) 
\prod_{1 \le i < j \le n} \frac{P^{k_j} - P^{k_i}}{P^j - P^i}.
$$
\end{theo}

\section{The combinatorial interpretation of $\alpha_{1,1}(n,0,S,f;k_1,\ldots,k_n)$}
\label{interpret}

If we specialize $P=1$, $Q=1$ and $m=0$ in Theorem~\ref{main} then we obtain
the following formula.
$$
\prod_{1 \le s < t \le n} \left( E_{k_t}  - {\Delta_{k_s}} {\Delta_{k_t}} 
\sum_{(i,j) \in S} f(i,j) {E_{k_t}}^i {E_{k_s}}^j \right) 
\prod_{1 \le i < j \le n} \frac{k_j - k_i}{j -i}
$$
By setting $s(X,Y) = - \sum\limits_{(i,j) \in S} f(i,j) \, Y^i X^j$, this is equal to \eqref{generalform} if 
$r(X,Y)$ is of the form given in \eqref{specialform}. Thus it suffices to interpret $\alpha_{1,1}(n,0,S,f;k_1,\ldots,k_n)$. The interpretation immediately follows from the recursion underlying this 
specialization, which reads as
\begin{multline*}
\sum_{(l_1,\ldots, l_{n-1})}^{(k_1,\ldots,k_n)} A(l_1,\ldots,l_{n-1})  
=   
\sum_{l_{n-1}=k_{n-1}}^{k_n} \sum_{(l_1,\ldots,l_{n-2})}^{(k_1,\ldots,k_{n-1})} 
A(l_1,\ldots,l_{n-1})    \\  + \, \sum_{(l_1,\ldots,l_{n-3})}^{(k_1,\ldots,k_{n-2})}
\sum_{(i,j) \in S} f(i,j) A(l_1,\ldots,l_{n-3},k_{n-1}+i,k_{n-1}+j).
\end{multline*} 

\medskip
 
For a finite subset $S \subseteq \mathbb{Z}^2$, we define an 
{\it $S$--triangle} to be a triangular array $(a_{i,j})_{1 \le i \le j \le n}$ (the entries are arranged in the same manner as those of monotone triangles) with the following 
properties: among the ``inner'' entries $(a_{i,j})_{1 \le i < j < n}$ of the
triangle, a certain set of non row adjacent
{\it special entries} is identified. If $a_{i,j}$ is such a special entry then $a_{i+1,j}$ and $a_{i+1,j+1}$ 
are said to be the parents of this special entry and we demand that there exists an $s=(r,t) \in S$ such that 
$(a_{i+1,j},a_{i+1,j+1})=(a_{i,j}+r,a_{i,j}+t)$. We say that the special entry
$a_{i,j}$ is associated to $s$. Observe that the second summand in the
recursion takes into account for these special entries. (In this case $k_{n-1}$ is the 
special entry.)
On the other hand and regarding the first summand in the recursion, for all $a_{i,j}$ with $i \not=1$ that are not a parent of a special entry, we demand that 
$a_{i-1,j-1} \le a_{i,j} \le a_{i-1,j}$ if $a_{i-1,j-1} \le a_{i-1,j}$ and $a_{i-1,j-1} > a_{i,j} > a_{i-1,j}$ 
if $a_{i-1,j-1} > a_{i-1,j}$; in the latter case we say that the pair $(a_{i-1,j-1},a_{i-1,j})$ is an {\it inversion}. We fix a function $f:S \to \mathbb{C}$ and define the weight of such an $S$--triangle as
$$
(-1)^\text{$\#$ of inversions} \prod_{s \in S} f(s)^\text{$\#$ of special entries associated to $s$}.
$$
With this definition, $\alpha_{1,1}(n,0,S,f;k_1,\ldots,k_n)$ is the sum of the weights of all $S$--triangles with 
bottom row $(k_1,\ldots,k_n)$.

\medskip

To give an example, observe that 
\begin{center}
\begin{tabular}{ccccccccccccccccc}
  &   &   &   &   &   &   &   & $6$ &   &   &   &   &   &   &   & \\
  &   &   &   &   &   &   & $5$ &   & $9$ &   &   &   &   &   &   & \\
  &   &   &   &   &   & $4$ &   & $10$ &   & $8$ &   &   &   &   &   & \\
  &   &   &   &   & $3$ &   & $4$ &   & $8^*$ &   & $9$ &   &   &   &   & \\
  &   &   &   & $7$ &   & $3^*$ &  &   $10$ &   & $7$   &  & $10$  &   &   &   & 
\end{tabular}
\end{center}
is a $\{(0,1),(2,0)\}$--triangle, where the special entries are marked with a star. The $S$--triangle has two inversions, one in the bottom row, i.e. $(10,7)$ and one in the middle row, i.e. $(10,8)$. Thus, if $f(0,1)=q_1$ and $f(2,0)=q_2$ then the weight of this $S$--triangle is $q_1 q_2$.

\medskip

Obviously, $\emptyset$--triangles with (weakly) increasing bottom row are simply Gelfand--Tsetlin patterns.
However, the notion of $\{(0,0)\}$--triangles does not coincide with the notion of monotone triangles. 
Still $\alpha_{1,1}(n,0,\{(0,0)\},-1;k_1,\ldots,k_n)$
is the number of monotone triangles with bottom row $(k_1,\ldots,k_n)$ if $k_1 < k_2 < \ldots < k_n$.
On the one hand and as mentioned above, this easily follows from the recursion underlying $\alpha_{1,1}(n,0,\{(0,0)\},-1;k_1,\ldots,k_n)$. 
On the other hand, this can also be shown  using the combinatorial interpretation that we have just introduced. For this purpose, fix 
$(k_1,\ldots,k_n) \in \mathbb{Z}^n$ with $k_1 <  k_2 < \dots < k_n$. Under this assumption, 
a $\{(0,0)\}$--triangle with bottom 
row $(k_1,\ldots,k_n)$ has no inversion. It suffices to show that the weighted sum $\Sigma$ with respect to 
$f \equiv -1$ over all 
$\{(0,0)\}$--triangles with bottom $(k_1,\ldots,k_n)$ that have at least one special entry is the negative 
of the number of all 
Gelfand--Tsetlin patterns with bottom row $(k_1,\ldots,k_n)$ that are not monotone triangles. If we ignore the marks of the special entries of a  $\{(0,0)\}$--triangle with bottom row $(k_1,\ldots,k_n)$ that appears in the sum $\Sigma$ then we clearly obtain a Gelfand--Tsetlin pattern with bottom row 
$(k_1,\ldots,k_n)$ that is not a monotone triangle. 
Thus, we have to show that for each fixed Gelfand--Tsetlin pattern that is not a monotone triangle the 
sum of weights of all $\{(0,0)\}$--triangles in $\Sigma$, whose unmarked version is equal to this fixed Geland--Tsetlin 
pattern is $-1$.  Indeed, suppose that $m$ is the number of row adjacent pairs that are equal in the fixed pattern. The candidates for the
special entries are those entries that are situated in a row below and in between such pairs. If we mark $k$ special entries in the Gelfand--Tsetlin pattern then the weight is 
$(-1)^k$ and, ignoring for a while the constraint that the special entries are
not suppossed to be row adjacent, there are clearly $\binom{m}{k}$
possiblities to do this. However, in this calculation the ``forbidden''
(i.e. row adjacent) markings cancel each other as $a_{i,j-1} = a_{i,j} = a_{i,j+1}$
implies $a_{i+1,j}=a_{i,j}=a_{i+1,j+1}$ and every forbidden marking that
contains $a_{i-1,j-1}$ and $a_{i-1,j}$ but not $a_{i,j}$ cancels with the
corresponding forbidden marking that includes $a_{i,j}$.
Thus, as $k$ ranges between $1$ and $m$, the multiplicity in question is 
$$
 \sum_{k=1}^{m} (-1)^k \binom{m}{k} =  \sum_{k=0}^{m} (-1)^k \binom{m}{k} - 1 = (1-1)^m - 1 = -1
$$
and the assertion follows. However, we want to emphasize that $\{(0,0)\}$--triangles (unlike 
monotone triangles) provide a combinatorial interpretation of $\alpha_{1,1}(n,0,\{(0,0)\},-1;k_1,\ldots,k_n)$ 
(and therefore of \eqref{formula}) if $k_i \ge k_{i+1}$ for an $i \in \{1,2,\ldots,n-1\}$.



\medskip

Although we think that it is of less interest, we conclude this section by discussing the case that $Q \not= 1$, i.e.
give a combinatorial interpretation of $\alpha_{1,Q}(n,0,S,f;k_1,\ldots,k_n)$.
 Here, we have different requirements for the $a_{i,j}$ with $i \not= 1$ that are not a parent of a special entry. Again we demand that $a_{i-1,j-1} \le a_{i,j} \le a_{i-1,j}$ if $a_{i-1,j-1} \le a_{i-1,j}$ and this entry contributes the weight $Q$ if  $a_{i-1,j-1} < a_{i,j} <  a_{i-1,j}$, the weight $1$ if $a_{i-1,j} = a_{i,j} < a_{i-1,j}$ or 
 $a_{i-1,j-1} < a_{i,j} = a_{i-1,j}$ and the weight $2-Q$ if 
$a_{i-1,j-1} = a_{i,j} =  a_{i-1,j}$. If   $a_{i-1,j-1} > a_{i-1,j}$ then
$a_{i-1,j-1} \ge a_{i,j} \ge a_{i-1,j}$ and this entry contributes the weight $-Q$ if 
$a_{i-1,j-1} > a_{i,j} > a_{i-1,j}$, otherwise it contributes $1 - Q$. In this case 
the total weight of a fixed $S$--triangle is the product of the $Q$--weights of its entries times 
$\prod\limits_{s \in S} f(s)^\text{$\#$ of special entries associated to $s$}$.

\section{The $Q$--enumeration of monotone triangles and the $P$--enumeration
  of weak monotone triangles}
\label{PQ} 

First we consider the $Q$--enumeration of monotone triangles. We claim that this weighted enumeration is obtained by specializing  
$m=0$, $P=1$, $S=\{(0,0)\}$ and $f \equiv -1$ in Theorem~\ref{main}. Indeed, in this case the recursion 
simplifies to 
\begin{multline*}
\sum_{(l_1,\ldots, l_{n-1})}^{(k_1,\ldots,k_n)} A(l_1,\ldots,l_{n-1})  \\
=  \left. \left( Q^{-1} (Q \id + \Delta_{k_n})  (Q E_{k^{*}_{n-1}} - \Delta_{k^{*}_{n-1}}) 
\sum_{l_{n-1}=k^{*}_{n-1}}^{k_n-1} \sum_{(l_1,\ldots,l_{n-2})}^{(k_1,\ldots,k_{n-1})} 
A(l_1,\ldots,l_{n-1})  \right) \right|_{k^{*}_{n-1} = k_{n-1}}  \\  - \, \sum_{(l_1,\ldots,l_{n-3})}^{(k_1,\ldots,k_{n-2})}
A(l_1,\ldots,l_{n-3},k_{n-1},k_{n-1}) \\
= \left. \left( (Q E_{k^{*}_{n-1}} + \Delta_{k_n} E_{k^{*}_{n-1}} - \Delta_{k^{*}_{n-1}} ) 
\sum_{l_{n-1}=k^{*}_{n-1}}^{k_n-1} \sum_{(l_1,\ldots,l_{n-2})}^{(k_1,\ldots,k_{n-1})} 
A(l_1,\ldots,l_{n-1})  \right) \right|_{k^{*}_{n-1} = k_{n-1}}  \\  - \, \sum_{(l_1,\ldots,l_{n-3})}^{(k_1,\ldots,k_{n-2})}
A(l_1,\ldots,l_{n-3},k_{n-1},k_{n-1}), 
\end{multline*}
which is furthermore equal to
\begin{multline*}
\sum_{(l_1,\ldots, l_{n-1})}^{(k_1,\ldots,k_n)} A(l_1,\ldots,l_{n-1})  \\
 =   Q \sum_{l_{n-1}=k_{n-1}+1}^{k_n-1} \sum_{(l_1,\ldots,l_{n-2})}^{(k_1,\ldots,k_{n-1})} 
A(l_1,\ldots,l_{n-1}) 
+ \sum_{(l_1,\ldots,l_{n-2})}^{(k_1,\ldots,k_{n-1})} 
A(l_1,\ldots,l_{n-2},k_{n})  \\
+ \sum_{(l_1,\ldots,l_{n-2})}^{(k_1,\ldots,k_{n-1})} 
A(l_1,\ldots,l_{n-2},k_{n-1})
 -  \sum_{(l_1,\ldots,l_{n-3})}^{(k_1,\ldots,k_{n-2})}
A(l_1,\ldots,l_{n-3},k_{n-1},k_{n-1}).
\end{multline*} 
It is straightforward to check that the specialization of the formula in Theorem~\ref{main} results in the formula of 
Theorem~\ref{Qweight}.

\medskip

Let us report on a subtility, which may be of importance when applying the ideas from \cite{fischernewproof} to 
evaluate $\alpha_{1,Q}(n,0,\{(0,0)\},-1;k_1,\ldots,k_n)$ at $(k_1,\ldots,k_n)=(1,2,\ldots,i-1,i+1,\ldots,n+1)$ 
in order to study a weighted refined enumeration of alternating sign matrices, see also Section~\ref{future}.
If we relax the condition 
of the strict increase of the rows of a monotone triangle to a possible weak increase in the bottom row then 
\eqref{formula} is the number of monotone triangles with bottom row $(k_1,\ldots,k_n)$ if 
$k_1 \le k_2 \le \ldots \le k_n$. However, this does not generalize to the $Q$--enumeration of monotone triangle as 
the recursion simplifies to 
\begin{multline*}   
\sum_{(l_1,\ldots, l_{n-1})}^{(k_1,\ldots,k_n)} A(l_1,\ldots,l_{n-1})  
= (2-Q)  \sum_{(l_1,\ldots,l_{n-2})}^{(k_1,\ldots,k_{n-1})} 
A(l_1,\ldots,l_{n-2},k_{n-1}) 
  \\  -  \sum_{(l_1,\ldots,l_{n-3})}^{(k_1,\ldots,k_{n-2})}
A(l_1,\ldots,l_{n-3},k_{n-1},k_{n-1})
\end{multline*} 
if $k_{n-1}=k_n$ and not to
\begin{multline*}   
\sum_{(l_1,\ldots, l_{n-1})}^{(k_1,\ldots,k_n)} A(l_1,\ldots,l_{n-1})  
=  \sum_{(l_1,\ldots,l_{n-2})}^{(k_1,\ldots,k_{n-1})} 
A(l_1,\ldots,l_{n-2},k_{n-1}) 
  \\  -  \sum_{(l_1,\ldots,l_{n-3})}^{(k_1,\ldots,k_{n-2})}
A(l_1,\ldots,l_{n-3},k_{n-1},k_{n-1}).
\end{multline*} 

\medskip 

Before we turn to the $P$--enumeration of weak monotone triangles, we will demonstrate how several $2$--enumeration of alternating--sign--matrix--structures follow from the $Q$--generating function. We start with the 
$2$--enumeration of ordinary $n \times n$ alternating sign matrices (with respect to the number of $-1$s in the alternating sign matrix), which was first established by Mills, Robbins and Rumsey in \cite{mills} as a corollary of their Theorem~2 (they have shown that it is given by 
$2^{\binom{n}{2}}$) long before the ordinary enumeration was finally
settled. Also our derivation shows that the $2$--enumeration is of a much simpler nature than the ordinary enumeration. Indeed, the generating function from Theorem~\ref{Qweight} simplifies to 
$$
\prod_{1 \le s < t \le n} \left( \id + E_{k_s} E_{k_t}  \right) 
\prod_{1 \le i < j \le n} \frac{k_j - k_i}{j - i}
$$
in this case.
Now, the crucial fact is that the operator $\prod\limits_{1 \le s < t \le n} \left( \id + E_{k_s} E_{k_t}  \right)$ is 
symmetric in $k_1,k_2,\ldots,k_n$ and, consequently, Lemma~1 from \cite{fischernewproof} implies that the generation function is 
equal to 
$
P(1,1,\ldots,1)  \prod\limits_{1 \le i < j \le n} \frac{k_j - k_i}{j - i},
$
where $P(X_1,X_2,\ldots,X_n) = \prod\limits_{1 \le s < t \le n} \left( \id + X_s X_t  \right)$ and, therefore, 
equal to 
\begin{equation}
\label{2}
2^{\binom{n}{2}}  \prod_{1 \le i < j \le n} \frac{k_j - k_i}{j - i}.
\end{equation}
(This result has already appeared implicitly in \cite{propp}.)
The $2$--enumeration is obtained by setting $k_i=i$ in this formula. If we specialize 
$(k_1,\ldots,k_{n-1})=(1,2,\ldots,l-1,l+1,\ldots,n)$ in the formula for the
$2$--enumeration of monotone triangles with bottom row $(k_{1},k_{2},\ldots,k_{n-1})$ then we obtain the $2$--enumeration of $n \times n$ alternating sign matrices where the unique $1$ in the first row is in the $l$--th column. That is
\begin{multline*}
2^{\binom{n-1}{2}}  \prod_{1 \le i < j \le n \atop i,j \not= l} (j - i) \prod_{1 \le i < j \le n-1} \frac{1}{j-i} 
= 2^{\binom{n-1}{2}} \prod_{1 \le i < j \le n} (j-i) \prod_{j=l+1}^{n} \frac{1}{j-l} \prod_{i=1}^{l-1} \frac{1}{l-i} 
\prod_{1 \le i < j \le n-1} \frac{1}{j-i} \\
= 2^{\binom{n-1}{2}} \prod_{i=1}^{n-1} (n-i) \frac{1}{(n-l)! (l-1)!} =  2^{\binom{n-1}{2}} \binom{n-1}{l-1}.
\end{multline*}
Similarly, one can reprove the $2$--enumeration of $(2n-1) \times (2n-1)$ vertically symmetric alternating sign matrices by specializing $(k_1,k_2,\ldots,k_{n})=(1,3,\ldots,2n-1)$ , which results in 
$2^{(n-1)(n-2))}$, see \cite{kuperberg}.

\medskip

In fact, it is possible to generalize the $2$--enumeration of $n \times n$ alternating sign matrices to certain partial alternating sign matrices. Let a 
partial $m \times n$ alternating sign matrix be an $m \times n$ matrix with entries $0, 1, -1$, where the entries 
$1$ and $-1$ alternate in each row and column, each row sum is $1$ and the first non--zero entry of each column is $1$ if there is any. (Such objects only exist if $m \le n$.) The well--known bijection between alternating sign matrices and monotone triangles shows that partial $m \times n$ alternating sign matrices are in bijection with monotone triangles 
with $m$ rows such that $1 \le k_1 < k_2 < \ldots < k_m \le n$. Thus and by \eqref{2}, the $2$--enumeration of these objects is given by 
\begin{equation}
\label{expr}
2^{\binom{m}{2}} \sum_{1 \le k_1 < k_2 < \ldots < k_m \le n}   \prod_{1 \le i < j \le m} \frac{k_j - k_i}{j - i} =
2^{\binom{m}{2}} \sum_{0 \le x_1 \le x_2 < \ldots \le x_m \le n-m}   \prod_{1 \le i < j \le m} \frac{x_j - x_i+j-i}{j - i}.
\end{equation}
As 
$
\prod\limits_{1 \le i < j \le m} \frac{x_j - x_i+j-i}{j - i}
$
is the number of semistandard tableaux of shape $(x_m,x_{m-1},\ldots,x_1)$, the sum on the right hand side of \eqref{expr} is the number of columnstrict plane partitions with at most $n-m$ columns and parts in $\{1,2,\ldots,m\}$. By the Bender--Knuth (ex--)Conjecture, see for instance \cite{fischermethod}, this number 
is equal to 
$
\prod\limits_{i=1}^{m} \frac{(n-m+i)_i}{(i)_i}
$
where $(a)_n = a(a+1) \cdots (a+n-1)$, and thus the $2$--enumeration of partial $m \times n$ alternating sign matrices is given by
$$
2^{\binom{m}{2}} \prod_{i=1}^{m} \frac{(n-m+i)_i}{(i)_i}.
$$

\medskip

For the $P$--enumeration of weak monotone triangles, the weighted enumeration is obtained by specializing  
$m=0$, $Q=1$, $S=\{(0,0)\}$ and $f \equiv -1$ in Theorem~\ref{main}, as the recursion simplifies to 
\begin{multline*}
\sum_{(l_1,\ldots, l_{n-1})}^{(k_1,\ldots,k_n)} A(l_1,\ldots,l_{n-1})  
=  \pe_{k_n} 
\sum_{l_{n-1}=k_{n-1}}^{k_n-1} P^{l_{n-1}} \sum_{(l_1,\ldots,l_{n-2})}^{(k_1,\ldots,k_{n-1})} 
A(l_1,\ldots,l_{n-1})    \\ 
- \, \sum_{(l_1,\ldots,l_{n-3})}^{(k_1,\ldots,k_{n-2})}
 A(l_1,\ldots,l_{n-3},k_{n-1},k_{n-1}) \\
= 
\sum_{l_{n-1}=k_{n-1}}^{k_n-1} P^{l_{n-1}} \sum_{(l_1,\ldots,l_{n-2})}^{(k_1,\ldots,k_{n-1})} 
A(l_1,\ldots,l_{n-1}) + 
 \sum_{(l_1,\ldots,l_{n-2})}^{(k_1,\ldots,k_{n-1})} 
A(l_1,\ldots,l_{n-2},k_{n})  \\ 
- \, \sum_{(l_1,\ldots,l_{n-3})}^{(k_1,\ldots,k_{n-2})}
 A(l_1,\ldots,l_{n-3},k_{n-1},k_{n-1})
\end{multline*} 
in this case. This easily implies that $\alpha_{P,1}(n,0,\{(0,0\},-1;k_1,\ldots,k_n)$ is the weighted enumeration of 
weak monotone triangles with bottom row $(k_1,\ldots,k_n)$ and with respect to the $P$--weight, which was defined 
in the introduction. Clearly, the formula in Theorem~\ref{main} simplifies to the formula in Theorem~\ref{Pweight} if 
we set  $m=0$, $Q=1$, $S=\{(0,0)\}$ and $f \equiv -1$.

\section{Proof of Theorem~\ref{main}}
\label{proof}

For fixed $S \subseteq \mathbb{Z}^2$ and $f:S \to \mathbb{C}$, we define the operator
$$
V_{x,y} = {\pqe_{x}} \pqid_y - {_ P\Delta_{x}} \, {\pd_{y}} \sum\limits_{(i,j) \in S} f(i,j) {\pe_{x}}^i {\pe_{y}}^{j}.
$$
The following lemma relates this operator to the summation operator defined in Section~\ref{master}.
\begin{lem} 
\label{key}
Suppose $B(l_1,\ldots,l_{n-1})$ is a function in $l_1,\ldots,l_{n-1}$ such that for all 
$i \in \{1,2,\ldots,n-2\}$ we have 
$V_{l_{i},l_{i+1}} B(l_1,\ldots,l_{n-1})=0$ if $l_{i}=l_{i+1}$. Then 
\begin{multline*}
\sum_{(l_1,\ldots,l_{n-1})}^{(k_1,\ldots,k_n)} {\pd_{l_1}} \dots {\pd_{l_{n-1}}} B(l_1,\ldots,l_{n-1}) \\
= 
\sum_{r=1}^{n} (-1)^{r-1} \prod_{s=1}^{r-1} \pqid_{k_s} \prod_{t=r+1}^{n} \pqe_{k_t}
B(k_1,\ldots,k_{r-1},k_{r+1},\ldots,k_n).
\end{multline*}
\end{lem}

{\it Proof.} We use induction with respect to $n$. For $n=0$ and $n=1$ there is nothing to prove. Suppose that $n \ge 2$. Then, by definition and by the 
induction hypothesis,  we have 
\begin{multline*}
\sum_{(l_1,\ldots,l_{n-1})}^{(k_1,\ldots,k_n)} {\pd_{l_1}} \cdots {\pd_{l_{n-1}}} B(l_1,\ldots,l_{n-1})  
= \left( Q^{-1} {\pqe_{k_n}} \pqid_{k^{*}_{n-1}} \sum_{l_{n-1}={k^{*}_{n-1}}}^{k_n-1}  P^{l_{n-1}} {\pd_{l_{n-1}}}  \right. \\  \sum_{r=1}^{n-1} (-1)^{r-1} 
\prod_{s=1}^{r-1} \pqid_{k_s} \prod_{t=r+1}^{n-1}  \pqe_{k_t}   \, 
 B(k_1,\ldots,k_{r-1},k_{r+1},\ldots,k_{n-1},l_{n-1})  \\
+ \sum_{r=1}^{n-2} (-1)^{r-1}   \prod_{s=1}^{r-1} \pqid_{k_s} \prod_{t=r+1}^{n-2}  \pqe_{k_t}
\sum_{(i,j) \in S} f(i,j) {\pe_{k_{n-1}}}^i  {\pe_{k^{*}_{n-1}}}^j \\ \left. \left.
{\pd_{k_{n-1}}} \, {\pd_{k^{*}_{n-1}}} B(k_1,\ldots,k_{r-1},k_{r+1},\ldots,k_{n-2},k_{n-1},k^{*}_{n-1}) \right) \right|_{k^{*}_{n-1}=k_{n-1}}.
\end{multline*}
This is obviously equal to 
\begin{multline*}
\sum_{r=1}^{n-1} (-1)^{r-1} \prod_{s=1}^{r-1} \pqid_{k_s} \prod_{t=r+1}^{n}  
 \pqe_{k_t}  \, 
 B(k_1,\ldots,k_{r-1},k_{r+1},\ldots,k_{n-1},k_{n}) \\
 + \sum_{r=1}^{n-1} (-1)^{r} \left. \left(  \prod_{s=1}^{r-1} \pqid_{k_s} \prod_{t=r+1}^{n-1}   \pqe_{k_t} \pqid_{k^*_{n-1}} 
 B(k_1,\ldots,k_{r-1},k_{r+1},\ldots,k_{n-1},k^*_{n-1}) \right) \right|_{k^*_{n-1}=k_{n-1}} \\
 + \sum_{r=1}^{n-2} (-1)^{r}  \prod_{s=1}^{r-1} \pqid_{k_s} \prod_{t=r+1}^{n-2}  \pqe_{k_t} 
 ((V_{k_{n-1},k^{*}_{n-1}} - \pqe_{k_{n-1}} \pqid_{k^{*}_{n-1}}) \\ 
\left.   B(k_1,\ldots,k_{r-1},k_{r+1},\ldots,k_{n-2},k_{n-1},k^{*}_{n-1}) ) \right|_{k^*_{n-1}=k_{n-1}}.
\end{multline*}
If we move the $(n-1)$--st summand of the second sum to the first sum and combine the third sum with the remainder of the second sum then we see that this is furthermore equal to \begin{multline*}
\sum_{r=1}^{n} (-1)^{r-1}  \prod_{s=1}^{r-1} \pqid_{k_s} \prod_{t=r+1}^{n}  
 \pqe_{k_t}
 B(k_1,\ldots,k_{r-1},k_{r+1},\ldots,k_{n-1},k_{n}) \\
 + \sum_{r=1}^{n-2} (-1)^r \left.  \prod_{s=1}^{r-1} \pqid_{k_s} \prod_{t=r+1}^{n-2}  \pqe_{k_t}
( V_{k_{n-1},k^*_{n-1}} \,  B(k_1,\ldots,k_{r-1},k_{r+1},\ldots,k_{n-2},k_{n-1},k^*_{n-1}) )
\right|_{k^*_{n-1}=k_{n-1}}
\end{multline*}
and the assertion follows. \qed

\medskip

In the following lemma we use the previous lemma to apply the summation operator to a special 
$B(l_1,\ldots,l_{n-1})$.

\begin{lem} 
\label{rec}
Let $m_2,m_3,\ldots,m_n \ge 0$ be integers and set $m_1=-1$. Then
\begin{multline*}
{ \sum_{(l_1,\ldots,l_{n-1})}^{(k_1,\ldots,k_n)}} \left( \prod_{1 \le s < t \le n-1} { V_{l_t,l_s}} \right)
\det_{1 \le i , j \le n-1} \pbinom{l_i}{m_{j+1}}  \\
= Q^{-n+1} \prod_{j=2}^{n} P^{m_j} \left( \prod_{1 \le s < t \le n} { V_{k_t,k_s}} \right)
\det_{1 \le i , j \le n} \pbinom{k_i}{m_j+1}.
\end{multline*}
\end{lem}

{\it Proof.} Let
$$
B(l_1,\ldots,l_{n-1}) = \left( \prod_{1 \le s < t \le n-1} { V_{l_t,l_s}} \right)
\det_{1 \le i , j \le n-1} \pbinom{l_i}{m_{j+1}+1}.
$$
First we show that $B(l_1,\ldots,l_{n-1})$ has the property, which is presumed in Lemma~\ref{key}. For that purpose
it suffices to show that
$$
(\id + S_{l_i,l_{i+1}}) { V_{l_i,l_{i+1}}} B(l_1,\ldots,l_{n-1}) = 0
$$
for all $i \in \{1,2,\ldots,n-2\}$, where $S_{x,y}$ is the {\it swapping operator}, defined as 
$S_{x,y} f(x,y)=f(y,x)$. This assertion follows, since, on the one hand,  
$$
{ V_{l_i,l_{i+1}}} \, { V_{l_{i+1},l_{i}}} 
\left( \prod_{1 \le s < t \le n-1 \atop (s,t) \not= (i,i+1)} { V_{l_t,l_s}} \right)
$$
is symmetric in $l_i$ and $l_{i+1}$ and therefore commutes with $S_{l_i,l_{i+1}}$ and, on the other hand,  
$\det\limits_{1 \le i , j \le n-1} \pbinom{l_i}{m_{j+1}+1}$ is antisymmetric in $l_i$ and $l_{i+1}$. Furthermore, by \eqref{diff},
we conclude that 
\begin{multline*}
{ \pd_{l_1}} \dots { \pd_{l_{n-1}}} B(l_1,\ldots,l_{n-1}) 
=  { \pd_{l_1}} \dots { \pd_{l_{n-1}}} \left( \prod_{1 \le s < t \le n-1} { V_{l_t,l_s}} \right)
\det_{1 \le i , j \le n-1} \pbinom{l_i}{m_{j+1}+1} \\
=  \left( \prod_{1 \le s < t \le n-1} { V_{l_t,l_s}} \right) 
{\pd_{l_1}} \dots {\pd_{l_{n-1}}} \det_{1 \le i , j \le n-1} \pbinom{l_i}{m_{j+1}+1} \\
=  \left( \prod_{1 \le s < t \le n-1} { V_{l_t,l_s}} \right) 
\det_{1 \le i , j \le n-1} { \pd_{l_i}} \pbinom{l_i}{m_{j+1}+1} 
=  \prod_{j=2}^{n} P^{-m_j} \left( \prod_{1 \le s < t \le n-1} { V_{l_t,l_s}} \right) 
\det_{1 \le i , j \le n-1}  \pbinom{l_i}{m_{j+1}}.
\end{multline*}
Therefore, by Lemma~\ref{key}, the left--hand side of the identity stated in the lemma 
is equal to 
$$
\prod_{j=2}^{n} P^{m_j} \sum_{r=1}^{n} (-1)^{r-1}  \prod_{s=1}^{r-1} \pqid_{k_s} \prod_{t=r+1}^{n} \pqe_{k_t}
B(k_1,\ldots,k_{r-1},k_{r+1},\ldots,k_n).
$$
By the definition of $B(l_1,\ldots,l_{n-1})$, this is equal to
\begin{multline*}
\prod_{j=2}^{n} P^{m_j} \sum_{r=1}^{n} (-1)^{r-1}
\left( \prod_{1 \le s < t \le n \atop s,t \not=r} { V_{k_t,k_s}} \right) 
  \prod_{s=1}^{r-1} \pqid_{k_s} \prod_{t=r+1}^{n} \pqe_{k_t} \\ 
\det_{1 \le i,j \le n-1} \left. \pbinom{l_i}{m_{j+1}+1} \right|_{(l_1,\ldots,l_{n-1})=(k_1,\ldots,\widehat{k_r},\ldots,k_n)}.
\end{multline*}
Since 
$$
{\pd_{k_r}} 
\det_{1 \le i,j \le n-1} \left. \pbinom{l_i}{m_{j+1}+1} \right|_{(l_1,\ldots,l_{n-1})=(k_1,\ldots,\widehat{k_r},\ldots,k_n)} =0,
$$
this is furthermore equal to
\begin{multline*}
Q^{-n+1} \prod_{j=2}^{n} P^{m_j}  \sum_{r=1}^{n} (-1)^{r-1}
\left( \prod_{1 \le s < t \le n \atop s,t \not=r} { V_{k_t,k_s}} \right) 
\prod_{s=1}^{r-1} { V_{k_{r},k_s}} \prod_{t=r+1}^{n} { V_{k_{t},k_{r}}} \\
 \det_{1 \le i,j \le n-1} \left. \pbinom{l_i}{m_{j+1}+1} \right|_{(l_1,\ldots,l_{n-1})=(k_1,\ldots,\widehat{k_r},\ldots,k_n)} \\
= Q^{-n+1} \prod_{j=2}^{n} P^{m_j}  \sum_{r=1}^{n} (-1)^{r-1}
\left( \prod_{1 \le s < t \le n} { V_{k_t,k_s}} \right) 
\det_{1 \le i,j \le n-1} \left. \pbinom{l_i}{m_{j+1}+1} \right|_{(l_1,\ldots,l_{n-1})=(k_1,\ldots,\widehat{k_r},\ldots,k_n)} \\
= Q^{-n+1} \prod_{j=2}^{n} P^{m_j} \left( \prod_{1 \le s < t \le n} { V_{k_t,k_s}} \right)  
\sum_{r=1}^{n} (-1)^{r-1}
\det_{1 \le i,j \le n-1} \left. \pbinom{l_i}{m_{j+1}+1} \right|_{(l_1,\ldots,l_{n-1})=(k_1,\ldots,\widehat{k_r},\ldots,k_n)}
\end{multline*}
This is now the right--hand side of the identity in the statement of the lemma, which is evident when 
expanding the determinant in the statement of the lemma with respect to the first column. \qed

\medskip

We define a quantity that is even more general than $\alpha_{P,Q}(n,m,S,f;k_1,\ldots,k_n)$. For $r \ge 1$ and 
$(m_1,\ldots,m_r) \in \mathbb{Z}_{\ge 0}^r$, let 
$$
\alpha_{P,Q}(1,(m_1,\ldots,m_r),S,f;k_1,\ldots,k_r) = 
\left( \prod_{1 \le s < t \le r} { V_{k_t,k_s}} \right)
\det_{1 \le i , j \le r} \pbinom{k_i}{m_j}
$$
and, for $n > 1$, let 
\begin{multline*}
\alpha_{P,Q}(n,(m_1,\ldots,m_r),S,f;k_1,\ldots,k_{n+r-1}) \\
= \sum_{(l_1,\ldots, l_{n+r-2})}^{(k_1,\ldots,k_{n+r-1})} \alpha_{P,Q}(n-1,(m_1,\ldots,m_r),S,f;l_1,\ldots,l_{n+r-2}).
\end{multline*}
By induction with respect to $n$, Lemma~\ref{rec} shows that $\alpha_{P,Q}(n,(m_1,\ldots,m_r),S,f;k_1,\ldots,k_{n+r-1})$
is equal to
\begin{multline*}
P^{\frac{1}{6} (n + 3r -3)(n-1) (n-2)  + (m_1 + m_2 + \ldots + m_r)(n-1)} Q^{-\binom{n}{2}}  \left( \prod_{1 \le s < t \le n+r-1} { V_{k_t,k_s}} \right) \\
\det_{1 \le i , j \le n+r-1} \pbinom{k_i}{[j < n](j-1)+[j \ge n](m_{j-n+1}+n-1)}.
\end{multline*}
The first statement in Theorem~\ref{main} is the special 
case $r=1$ and $m_1=m$. The second statement follows as 
$$
\det_{1 \le i, j \le n} \pbinom{k_i}{j-1} = P^{\binom{n}{2}} \prod_{1 \le i < j \le n} 
\frac{P^{k_j} - P^{k_i}}{P^j - P^i}
$$
by the $q$--Vandermonde determinant evaluation.

\section{Remarks and further projects}
\label{future}

The starting point for this paper was \cite{fischermonotone}, where we have studied the
recursion underlying the counting function for monotone triangles with prescribed bottom row. In the present paper, we have considered a generalized recursion, which we have obtained by carefully introducing 
a number of new parameter, namely $m, P, Q$, a finite subset $S \subseteq \mathbb{Z}^2$ and a 
function $f:S \to \mathbb{C}$, in the original recursion. The analysis of this generalized 
recursion was possible as we have finally noticed that the analysis of the original recursion can be simplified significantly.

\medskip

With the exception of $m$, all new parameters have been used to either obtain weighted enumerations of monotone 
triangles, respectively weak monotone triangles or a combinatorial interpretation of a generalization of 
\eqref{formula}. However, the parameter $m$ is also of special interest since it offers the possibility to ``control'' top 
and bottom row of a monotone triangle -- so far we were only able to ``control'' either row. Indeed, for fixed 
$n \ge 1$ and $1 \le i \le n$, let $(c_{p,q})_{p,q \ge 0}$ be complex coefficients, where almost 
all coefficients are zero, 
such that the polynomial $\sum\limits_{p,q \ge 0} c_{p,q} \binom{k-p}{q}$ (in $k$) vanishes for all 
$k \in \{1,2,\ldots, n\} \setminus \{i\}$ and is equal to $1$ for $k=i$. Then, it is not hard to see that
$$
\sum_{p,q \ge 0} c_{p,q} \, \alpha_{1,1}(n,q,\{(0,0)\},-1;k_1-p,\ldots,k_n-p)
$$
is the number of monotone triangles with bottom row $(k_1,\ldots,k_n)$ and top row $i$. (Here, we need the 
fact that 
$$
\sum_{(l_1,\ldots,l_{n-1})}^{(k_1,\ldots,k_n)} A(l_1-p,l_2-p,\ldots,l_{n-1}-p) =
\sum_{(l_1,\ldots,l_{n-1})}^{(k_1-p,\ldots,k_n-p)} A(l_1,l_2,\ldots,l_{n-1}).)
$$
In a 
forthcoming paper, we want to use the methods from \cite{fischernewproof} to
attack a certain doubly refined 
enumeration of alternating sign 
matrices, namely study the number of $n \times n$ alternating sign matrices where the unique $1$ in the top 
row is in column $i$ and the unique $1$ in the bottom row is in column
$j$. (Cleary, this number is equal to 
$$
\sum_{p,q \ge 0} c_{p,q} \, \alpha_{1,1}(n-1,q,\{(0,0)\},-1;1-p,2-p,\ldots,j-1-p,j+1-p,j+2-p,\ldots,n-p).)
$$
Note that this doubly refined enumeration of alternating sign matrices has already been
considered in \cite{stroganov}. Promising computerexperiments show that 
$\alpha_{1,1}(n-1,q,\{(0,0)\},-1;1-p,\ldots,j-1-p,j+1-p,\ldots,n-p)$ is in fact ``round'' 
(i.e. has only relatively small prime factors) for certain choices of $p$ and $q$. For instance, we have 
worked out the following conjecture for $q=1$, 
\begin{multline*}
\alpha_{1,1}(n-1,1,\{(0,0)\},-1;1-p,2-p,\ldots,j-1-p,j+1-p,j+2-p,\ldots,n-p) \\
= (j-p) A_{n-1,j} + A_n \frac{(n-j+1)_{2j-3} n (n-2j+1) (n+j-1)}{(2n - j - 1) (2n - j + 1)_{j-1} (j-1)!},
\end{multline*}
where $A_n = \prod\limits_{j=0}^{n-1} \frac{(3j+1)!}{(n+j)!}$ is the number of $n \times n$ alternating sign matrices 
and $$A_{n,i} = \binom{n+i-2}{i-1} \frac{(2n-i-1)!}{(n-i)!} \prod_{j=0}^{n-2} \frac{(3j+1)!}{(n+j)!}$$ is the 
number of $n \times n$ alternating sign matrices that have a $1$ which is situated in the first row and 
$i$--th column.

\medskip

Of course, it is also of interest to  apply the ideas from 
\cite{fischernewproof} to obtain informations on the special evaluations at 
$(k_1,\ldots,k_n)=(1,2,\ldots,i-1,i+1,\ldots,n+1)$ 
of the other generalizations of \eqref{formula}. For instance, the evaluation of the $Q$--enumeration 
in Theorem~\ref{Qweight} at $(k_1,\ldots,k_n)=(1,2,\ldots,n)$ is the weighted enumeration of $n \times n$ 
alternating sign matrices with respect to the number of $-1$s in the alternating sign matrix and, similary, 
the evaluation at  $(k_1,\ldots,k_n)=(1,3,\ldots,2n-1)$ is the weighted enumeration of 
$(2n-1) \times (2n-1)$ vertically symmetric alternating sign matrices. So far, there do not exist formulas 
for these generating functions (it is likely that it is rather difficult to come up with a simple formula as 
the generating functions do not seem to be 
``round''), however Kuperberg \cite[Theorem~4]{kuperberg} proved a number of factorizations regarding 
generating functions of this type. Is it possible to use the methods from \cite{fischernewproof} to deduce 
refined versions of these relations?

\medskip

Another natural question is whether the new insight in the recursion underlying \eqref{formula} leads to further generalizations of the formula. A project along these lines will be the following: for $1 \le i \le n$, let 
a monotone $(i,n)$--trapezoid be a monotone triangle with the first $i-1$ rows removed. We want to study the number of monotone $(i,n)$--trapezoids with prescribed bottom row $(k_1,\ldots,k_n)$. Obviously, the underlying 
recursion is the same as those for monotone triangles with prescribed bottom row. The difference (and the 
difficulty) lies in the inital condition. Secondly, we want to remark 
that our $P$--enumeration of weak monotone triangles with prescribed bottom row is a result of our efforts 
to obtain a weighted enumeration of monotone triangles, where the weight of a
given monotone triangle is equal (or at least related) to $P$ raised to the power of the sum of entries. 
Is it possible to write down a closed (operator) formula for a
generating function of this type, or is the $P$--enumeration of weak monotone
triangles already the best we can achieve in this respect? Thirdly, it should
be mentionend that the result regarding the combinatorial interpretation of
\eqref{generalform} emerged during our efforts to find other objects whose
counting functions can be expressed by an operator formula. The (somehow
inverse) strategy, which finally led to the result, was to start with an
operator formula and to search for objects whose weighted enumeration is given
by the formula. The solution given in Section~\ref{interpret} is in the sense
not satisfactory as it does not lead to a ``plain'' enumeration when
specializing the weights and, moreover, the definition of $S$--triangles is a
bit involved. Clearly, it would be of interest to search for other operators
formulas and (simpler) combinatorial objects that are enumerated by these formulas by varying the operators and/or the (factorizing) polynomial to which the operators are applied.

\bigskip \noindent
\textsc{
\!\!Ilse Fischer\\
Institut f\"ur Mathematik, Universit\"at Klagenfurt \\
9020 Klagenfurt, and \\
Fakult\"at f\"ur Mathematik, Universit\"at Wien \\
1090 Wien, Austria \\
}
\texttt{Ilse.Fischer@univie.ac.at}

\end{document}